\documentclass[12pt]{amsart}
\usepackage{amsmath}
\usepackage{amssymb, amsfonts, amsthm}
\usepackage{mathtools}
\usepackage{amssymb}
\usepackage{ifthen}
\usepackage{graphicx}
\usepackage{float}
\usepackage{caption}
\usepackage{subcaption}
\usepackage{cite}
\usepackage{mathrsfs}
\usepackage{amsfonts}
\usepackage{amscd}
\usepackage{amsxtra}
\usepackage{color}
\usepackage{bm}
\numberwithin{equation}{section}

\newtheorem{theorem}{Theorem}[section]
\newtheorem{lemma}[theorem]{Lemma}
\newtheorem{corollary}[theorem]{Corollary}

\theoremstyle{definition}

\newtheorem{examples}[theorem]{Examples}

\newcounter{alphabet}
\newcounter{tmp}

\def\be{\begin{equation}}
\def\ee{\end{equation}}

\newcommand{\ben}{\begin{enumerate}}
\newcommand{\een}{\end{enumerate}}

\newcommand{\blem}{\begin{lemma}}
\newcommand{\elem}{\end{lemma}}
\newcommand{\bthm}{\begin{theorem}}
\newcommand{\ethm}{\end{theorem}}
\newcommand{\bcor}{\begin{corollary}}
\newcommand{\ecor}{\end{corollary}}
\newcommand{\beg}{\begin{exam}}
\newcommand{\eeg}{\end{exam}}
\newcommand{\begs}{\begin{examples}}
\newcommand{\eegs}{\end{examples}}
\newcommand{\bdefe}{\begin{defn}}
\newcommand{\edefe}{\end{defn}}
\newcommand{\bprob}{\begin{prob}}
\newcommand{\eprob}{\end{prob}}
\newcommand{\bques}{\begin{ques}}
\newcommand{\eques}{\end{ques}}
\newcommand{\bei}{\begin{itemize}}
\newcommand{\eei}{\end{itemize}}
\newcommand{\bcon}{\begin{conj}}
\newcommand{\econ}{\end{conj}}
\newcommand{\bop}{\begin{op}}
\newcommand{\eop}{\end{op}}

\newcommand{\bas}{\begin{assertion}}
\newcommand{\eas}{\end{assertion}}

\newcommand{\bfa}{\begin{fact}}
\newcommand{\efa}{\end{fact}}

\newcommand{\bca}{\begin{ca}}
\newcommand{\eca}{\end{ca}}

\newcommand{\bst}{\begin{step}}
\newcommand{\est}{\end{step}}

\newcommand{\bsca}{\begin{sca}}
\newcommand{\esca}{\end{sca}}

\newcommand{\bcl}{\begin{cl}}
\newcommand{\ecl}{\end{cl}}

\newcommand{\bmlem}{\begin{mlem}}
\newcommand{\emlem}{\end{mlem}}

\newcommand{\bscl}{\begin{scl}}
\newcommand{\escl}{\end{scl}}

\newcommand{\bcons}{\begin{conjs}}
\newcommand{\econs}{\end{conjs}}

\newcommand{\bprop}{\begin{prop}}
\newcommand{\eprop}{\end{prop}}

\newcommand{\br}{\begin{remark}}
\newcommand{\er}{\end{remark}}
\newcommand{\brs}{\begin{rems}}
\newcommand{\ers}{\end{rems}}
\newcommand{\bo}{\begin{obser}}
\newcommand{\eo}{\end{obser}}
\newcommand{\bos}{\begin{obsers}}
\newcommand{\eos}{\end{obsers}}
\newcommand{\bpf}{\begin{pf}}
\newcommand{\epf}{\end{pf}}
\newcommand{\ba}{\begin{array}}
\newcommand{\ea}{\end{array}}
\newcommand{\beq}{\begin{eqnarray}}
\newcommand{\beqq}{\begin{eqnarray*}}
\newcommand{\eeq}{\end{eqnarray}}
\newcommand{\eeqq}{\end{eqnarray*}}

\newcommand{\ra}{\to}

\newcounter{minutes}
\divide\time by 60
\newcounter{hours}
\multiply\time by 60 \addtocounter{minutes}{-\time}

\begin{document}
	\title[Bohr-type inequalities for holomorphic functions]
	{The Bohr-type inequalities for holomorphic functions with lacunary series in complex Banach space}

\thanks{
File:~\jobname .tex,
          printed: \number\day-\number\month-\number\year,
          \thehours.\ifnum\theminutes<10{0}\fi\theminutes}

	\author[S. Kumar]{Shankey Kumar}
	\address{Shankey Kumar, Department of Mathematics, Indian Institute of Technology Madras, Chennai, 600036, India.}
	\email{shankeygarg93@gmail.com}
	
	\author[S. Ponnusamy]{Saminathan Ponnusamy}
	\address{Saminathan Ponnusamy, Department of Mathematics, Indian Institute of Technology Madras, Chennai, 600036, India.}
\address{Lomonosov Moscow State University, Moscow Center of Fundamental and Applied Mathematics, Moscow, Russia.}
	\email{samy@iitm.ac.in}

\author[G. B. Williams]{G. Brock Williams}
	\address{G. Brock Williams, Department of Mathematics and Statistics, Texas Tech University, Lubbock, TX, 79409, USA.}
\email{Brock.Williams@ttu.edu}

	
\begin{abstract}
	In this paper, we study the Bohr inequality with lacunary series to the single valued (resp. vector-valued) holomorphic function defined in unit ball of finite dimensional Banach sequence space. Also, we extend the Bohr inequality with an alternating series to the higher-dimensional space.
\end{abstract}
	
\maketitle
\pagestyle{myheadings}
\markboth{S. Kumar, S. Ponnusamy and G. B. Williams}{Bohr-type inequalities for holomorphic functions}

\section{Introduction and Preliminaries}\label{LLP-sec1}

\subsection{The classical Bohr inequality for the class $\mathcal{B}$} We denote by $\mathbb{D}:=\{z\in \mathbb{C}:|z|<1 \}$ the unit disk in the complex plane.
Let $\mathcal{H}$ be the class of all holomorphic functions defined on $\mathbb{D}$. We set $\mathcal{B}=\{f\in\mathcal{H}: |f(z)|\leq 1\}$.
Let us first recall a remarkable result of Bohr \cite{B1914} that opens up a new avenue for  research in geometric function theory.

\vspace{8pt}
\noindent
{\bf Theorem A.}\quad
{\it
Let $f\in \mathcal{B}$ be of the form $f(z)=\sum_{k=0}^{\infty} a_{k} z^{k}$. Then the Bohr sum $B_f(r)$
satisfies the inequality
$$
B_f(r):=\sum_{k=0}^{\infty}|a_k|r^k\le1\quad for \quad |z|=r\le\frac{1}{3}
$$
and the constant $1/3$ is best possible.}

\vspace{8pt}

Originally, Bohr obtained this inequality only for $r\le 1/6$, but later M. Riesz, I. Schur and F. W. Wiener independently proved that it holds in this form, and the number $1/3$ is called the Bohr radius. Several other proofs are also known in the literature. Note that there is no extremal function such that the Bohr radius is precisely $1/3$. See \cite[Corollary 8.26]{GarMasRoss-2018}. The Bohr inequality contemplates many generalizations and
applications.  In 1995, Dixon \cite{Dixon-95-7} used the Bohr inequality in connection with the
long-standing open problem of characterizing Banach algebras satisfying the von Neumann inequality.
In addition, Theorem~A
was extended to alternating series $A_f(r)=\sum_{k=0}^{\infty}(-1)^k|a_k|r^k$, by Ali et al. (cf. \cite{ABS2017}).
More precisely they have shown that $|A_f(r)|\leq1$ for $r\le1/\sqrt{3}$ under the assumption of Theorem A.
Another  natural question was to discuss the asymptotic  behaviour of the Bohr sum $B_f(r)$. This has led to the search for the best constant $C(r)\geq 1$ such that $B_f(r)\leq C(r)$. Bombieri \cite{Bomb-1962} showed that if $f\in \mathcal{B}$, then
$$
B_f(r)\leq \frac{3-\sqrt{8(1-r^2)}}{r} ~\mbox{ for }~ \frac{1}{3}\leq r\leq \frac{1}{\sqrt{2}}.
$$
For an alternate proof of this inequality, we refer to the recent paper \cite{KayPon_AAA18}.
Moreover, the question raised by Djakov and Ramanujan \cite{DjaRaman-2000} about $p$-Bohr radius was answered affirmatively in \cite{KayPon_AAA18}.
However, in the year 2004, Bombieri and Bourgain  \cite{BombBour-2004} established that $B_f(r)< \frac{1}{\sqrt{1-r^{2}}}$ holds for $r> 1/\sqrt{2}$
which in turn implies that $C(r) \asymp (1-r^{2})^{-1/2}$ as $r\to 1$. In the same article, the authors proved that for a given $\varepsilon >0$ there
exists a constant $c$ depending on $\varepsilon$ such that
$$
B_f(r) \geq (1-r^2)^{-1/2} - \left( c \log \frac{1}{1-r} \right)^{3/2+\varepsilon} ~\mbox{ as $r\rightarrow 1$}.
$$
Some recent results on this topic including refinements and generalizations  may be found from \cite{V2022,KKP2017,KayPon_AAA18,LLP2021P,LSX2018,PVW2020,PVW2022,PVW2021}.

\subsection{Multi-dimensional Bohr's inequality}
In the recent years, many authors paid attention to multidimensional generalizations of Bohr's theorem and drew many conclusions.
For example, denote an $n$-variables power series by $\sum\limits_{\alpha}a_{\alpha}z^{\alpha}$ with the standard multi-index notation; $\alpha=(\alpha_{1}, \alpha_{2},\ldots,\alpha_{n})$, where $\alpha_{j}\in
\mathbb{N}_0:=\mathbb{N}\cup\{0\}$,  $\mathbb{N}:=\{1,2,\ldots \}$ ($1\leq j\leq n$), $|\alpha|$ denotes the sum $\alpha_{1}+\alpha_{2}+\cdots+\alpha_{n}$ of its components,
$\alpha!=\alpha_{1}! \alpha_{2}! \cdots\alpha_{n}!$, $z=(z_{1},z_{2},\ldots ,z_{n})\in \mathbb{C}^n$, and $z^{\alpha}=z^{\alpha_{1}}z^{\alpha_{2}} \cdots z^{\alpha_{n}}$. The $n$-dimensional Bohr radius $K_{n}$ is the largest number
such that if $\sum\limits_{\alpha}a_{\alpha}z^{\alpha}$ converges in the $n$-dimensional unit polydisk
$\mathbb{D}^n$
such that $\Big|\sum\limits_{\alpha}a_{\alpha}z^{\alpha}\Big|<1$ in $\mathbb{D}^n$, then $\sum\limits_{\alpha}|a_{\alpha}z^{\alpha}| \leq 1$ for
$\operatorname{max}_{1\leq j\leq k}|z_{j}|\leq K_{n}.$ In 1997, Boas and Khavinson \cite{BK1997} showed that for $n>1$, the $n-$dimensional Bohr radius $K_n$ satisfies
$$ \frac{1}{3\sqrt{n}}< K_{n}< 2 \sqrt{ \frac{\log n}{n}}.
$$
This article became a source of inspiration for many subsequent investigations including connecting the asymptotic behaviour of $K_n$ to problems in the geometry of Banach
spaces (cf. \cite{DGMP19}).  However determining  the exact value of the Bohr radius $K_{n}$, $n>1$, remains an open problem. In 2006, Defant and Frerick \cite{DF}
improved the lower bound as $K_n \geq c \sqrt{\log n /(n \log \log n)}$ whereas Defant et al. \cite{DFOOS} used the hypercontractivity of the
polynomial Bohnenblust-Hille inequality and showed that
$$K_n = b_n \sqrt{ \frac{\log n}{n}} ~\mbox{ with }~ \frac{1}{\sqrt{2}} +o(1)\leq b_n\leq 2.
$$
In 2014, Bayart et al. \cite{Bay} established the   asymptotic behaviour of $K_n$ by showing that
$$\lim_{n\ra \infty} \frac{K_n}{\sqrt{ \frac{\log n}{n}} }= 1.
$$
We would like to mention that
Djakov and Ramanujan \cite{DjaRaman-2000}, and Blasco \cite{B2009} have studied the asymptotic behavior of the holomorphic functions with $p$-norm as $r\rightarrow1$ in $\mathbb{D}^n$
and Banach spaces.
Aizenberg \cite{A1998,A2005} mainly generalized Carath\'eodory's inequality for functions holomorphic in $\mathbb{C}^n$. In 2021, Liu and Ponnusamy \cite{LS2021}
have established several multidimensional analogues of refined Bohr's inequality for holomorphic functions on complete circular domain in $\mathbb{C}^n$.
Other aspects and promotion of Bohr inequality in higher dimensions can be obtained from \cite{BDK5,DePre06,GK2003,HHK2009,LW2007}. Moreover, research on Dirichlet series in higher
dimensions is also very popular recently (see \cite{DGMP19}).

\subsection{Generalizations and Refinements of Bohr's inequality for the disk}
Recently, Kayumov and Ponnusamy \cite{KKP2017}, and Ponnusamy et al. \cite{PVW2020} established several refined
versions and improved versions of Bohr's inequality in the planer case. See also \cite{PVW2020,PVW2021,PVW2022}.


\vspace{8pt}
\noindent
{\bf Theorem B.}(\cite{LS2021,PVW2020})\quad
{\it
	For $f\in \mathcal{B}$, and $f(z)=\sum_{k=0}^{\infty} a_{k} z^{k}$, we have
	$$
	\sum_{k=1}^{\infty}|a_k|r^k+\left(\frac{1}{1+|a_0|}+\frac{r}{1-r}\right)\sum_{k=1}^{\infty}|a_k|^2r^{2k} \leq \frac{r}{1-r}(1-|a_0|^2).
	$$
}

\vspace{8pt}

In the planar case, generalizations of Theorem B are abundant, but they are still limited in the higher-dimensional space. In fact, the alternating series version of this result is contained in the following result which is indeed itself a special case of  \cite[Theorem 5(I)]{LLP2021P}.

\vspace{8pt}
\noindent
{\bf Theorem C.}\quad
{\it Suppose that $m\in\mathbb{N}_0$, $p\in\mathbb{N}$ and $0\le m\le p$. Let $f\in \mathcal{B}$, and
 $f(z)=\sum_{k=0}^{\infty} a_{kp+m} z^{kp+m}$. If 	$p$ is odd, then
	$$
	\left |\sum_{k=1}^{\infty}(-1)^{kp+m}|a_{kp+m}|r^{kp+m} +(-1)^{m+p} \frac{r^{p+m}}{1-r^{2p}}\sum_{k=0}^{\infty}|a_{kp+m}|^2r^{2kp}\right | \le 1
	$$
	holds for $|z|=r\le r_*$, where $r_{*}$ is the unique root in $(0, 1)$ of $r^{p}(r^{p}+r^{m})-1=0$. This result is sharp.
}

\vspace{8pt}

\subsection{Bohr radius in higher dimensional setting}
In 2019, Liu and Liu \cite{LL2020} used the Fr\'echet derivative to establish the Bohr inequality of norm type for holomorphic mappings with lacunary series on the unit polydisk in $\mathbb{C}^n$ under some restricted conditions. 
The relevant properties of the Fr\'echet derivative can be seen below (cf. \cite{GK2003}).

Let $X$ and $Y$ be two complex Banach spaces with respect to the norms $\Vert\cdot\Vert_1$ and $\Vert\cdot\Vert_2$, respectively. For convenience, we denote both norms by $\Vert\cdot\Vert$, when there is no possibility of confusion with the underlying Banach spaces. We set $\mathbb{B}_X:=\{x\in X:\|x\|<1\}$. Let $\Omega^*$ be a domain in $X$, and let $H(\Omega^*,Y)$ denote the set of all holomorphic mappings from $\Omega^*$ into $Y$. It is well-known (cf. \cite{GK2003}) that if $f\in H(\Omega^*,Y)$, then
\begin{equation*}
f(x)=\sum_{k=0}^{\infty} \frac{1}{k!}D^kf(0)(x^k)
\end{equation*}
for all $x$ in some neighborhood of $0\in \Omega^*$, where $D^kf(x)$ is the $k^{th}$-Fr\'echet derivative of $f$ at $x$, and for $k\ge1$, we have
\begin{equation*}
D^kf(0)(x^k)=D^kf(0)(\underbrace{x,x,\dots,x}_{k}).
\end{equation*}
Moreover, if $k=0$, then $D^0f(0)(x^0)=f(0)$.

\section{Key lemmas and their Proofs}\label{LLP-sec2}
In order to establish our main results, we need the following lemmas which play a key role in proving the subsequent results in Section \ref{LLP-sec3}.
The above theorem has been generalized in \cite[Lemma 4 ]{LLP2021P}.

\vspace{8pt}
\noindent
{\bf Lemma D.}\quad
{\it Suppose that $m\in\mathbb{N}_0$, $p\in\mathbb{N}$ and $0\le m\le p$.
If $f\in \mathcal{B}$ and $f(z)=\sum_{k=0}^{\infty}a_{pk+m}z^{pk+m}$, then we have
\begin{equation} \label{LemD*}
\sum_{k=1}^{\infty}|a_{(2k-1)p+m}| r^{(2k-1)p} + \frac{r^{2p}}{1-r^{2p}}\sum_{k=0}^{\infty}|a_{kp+m}|^2 r^{(2k-1)p}\leq\frac{r^p}{1-r^{2p}}
\end{equation}
for $r\in [0, 1)$. This result is sharp. 
Moreover,
\begin{eqnarray}\label{liu21b}
\sum_{k=1}^{\infty}|a_{2kp+m}| r^{2kp} + \left(\frac{1}{1+|a_m|}+\frac{r^{2p}}{1-r^{2p}}\right)\sum_{k=1}^{\infty}|a_{kp+m}|^2 r^{2kp}
 \leq (1-|a_m|^2)\frac{r^{2p}}{1-r^{2p}}
\end{eqnarray}
holds for $r\in [0, 1)$. This result is sharp for $f(z)=z^m\left(\frac{a-z^p}{1-az^p}\right)$ with $a\in[0,1)$. %
}

\vspace{8pt}

Let $n\in \mathbb{N}$, $t\in[1,\infty)$, and $B_{\ell_t^n}$ be the set defined as the collection of complex vectors
$z=(z_1,z_2,\dots,z_n)\in \mathbb{C}^n$ satisfying $\sum_{i=1}^{n}|z_i|^t<1$. This set constitutes the open unit ball in the complex Banach space $\ell_t^n$
where the norm $\|z\|_t$ of
$z$ is given by $\big (\sum_{i=1}^{n}|z_i|^t\big )^{1/t}$. In the special case of $B_{\ell_\infty^n}$,
the set represents the unit polydisk in $\mathbb{C}^n$ denoted as $B_{\ell_\infty^n}:=\mathbb{D}^n$, where $|z_i|<1$ for $1\leq i\leq n$. The norm of $z\in \ell_\infty^n$ is defined as $\|z\|:=\max \{|z_i|:1\leq i\leq n\}$.
Note that the unit disk $\mathbb{D}$ is equivalent to $B_{\ell_t^1}$.

\blem\label{LLP-lem2}
Suppose that $m,\, p\in\mathbb{N},$ $1\le m\le p$, $1\leq t \leq \infty$, $f  \in H(B_{\ell_t^n},\overline{\mathbb{D}}^n)$ and
$$
f(z)=zg(z)=\frac{D^mf(0)(z^m)}{m!}+\sum_{k=1}^{\infty}\frac{D^{kp+m}f(0)(z^{kp+m})}{(kp+m)!},
$$
where $g\in H(B_{\ell_t^n},\mathbb{C})$.
Then
\begin{equation*}
\sum_{k=1}^{\infty} \frac{\left \Vert D^{2kp+m}f(0)(z^{2kp+m})\right \Vert_\infty}{(2kp+m)!}
\le\frac{\Vert z\Vert_t ^{2p-m}}{1-\Vert z\Vert_t ^{2p}}\left [ \Vert z\Vert_t^{2m}-\left(\frac{\Vert D^m f(0)(z^m)\Vert_\infty}{m!}\right)^2\right ]
\end{equation*}
holds for all $\Vert z\Vert_t=r\in[0,1)$.
\elem

\noindent
{\bf Proof.}\quad The proof of this lemma follows if we adopt the same approach as \cite[Lemma 2.1]{LLP23}.
To do this, we fix $z\in B_{\ell_t^n}\backslash\{ 0 \}$, and let $z_0=\frac{z}{\Vert z \Vert_t}$. Then $z_0\in\partial B_{\ell_t^n}$. Define
$j$ such that $|z_j|=\Vert z\Vert_\infty=\max_{1\le l\le n}\{ |z_l| \}.$ Because $f(z)=zg(z)$, a simple calculation yields
\begin{equation*}
\frac{D^kf_j(0)(z_0^k)}{k!}=\frac{D^{k-1}g(0)(z_0^{k-1})}{(k-1)!}\frac{z_j}{\Vert z\Vert_t} ~\mbox{ for each  $k\in\mathbb{N}$.}
\end{equation*}

For $\lambda\in \mathbb{D}$, let $h_j(\lambda)=f_j(\lambda z_0)$. Then $h_j\in H(\mathbb{D},\overline{\mathbb{D}})$ and
\begin{eqnarray*}
h_j(\lambda)&=&\lambda^m\left[\frac{D^mf_j(0)(z_0^m)}{m!}+\sum_{k=1}^{\infty}\frac{D^{kp+m}f_j(0)(z_0^{kp+m})}{(kp+m)!}\lambda^{kp}\right]\\
&=&\lambda^m\left[\frac{D^{m-1}g(0)(z_0^{m-1})}{(m-1)!}\frac{z_j}{\Vert z\Vert_t}+\sum_{k=1}^{\infty}\frac{D^{kp+m-1}g(0)(z_0^{kp+m-1})}{(kp+m-1)!}\frac{z_j}{\Vert z\Vert_t}\lambda^{kp}\right]\\
&=:&\lambda^m\varphi(\lambda^p), 
\end{eqnarray*}
where $\varphi(t)=b_0+\sum_{k=1}^{\infty}b_k t^k$ with
$$b_0=\frac{D^{m-1}g(0)(z_0^{m-1})}{(m-1)!}\frac{z_j}{\Vert z\Vert_t} \ \mbox{ and } \ b_k=\frac{D^{kp+m-1}g(0)(z_0^{kp+m-1})}{(kp+m-1)!}\frac{z_j}{\Vert z\Vert_t}.$$

Note that $\varphi \in H(\mathbb{D},\overline{\mathbb{D}})$ and so, $|b_0|\le1$. By Wiener's inequality, it follows that $|b_k|\le1-|b_0|^2$ for all $k\ge1$.   As $|z_j|=\Vert z\Vert_\infty=\max_{1\le l\le n}\{|z_l|\}$, we have
\begin{equation*}
 \frac{|D^mf_l(0)(z_0^m)|}{m!}=\frac{|D^{m-1}g(0)(z_0^{m-1})|}{(m-1)!}\frac{|z_l|}{\Vert z\Vert_t}\le\frac{|D^{m-1}g(0)(z_0^{m-1})|}{(m-1)!}\frac{|z_j|}{\Vert z\Vert_t}\le 1
\end{equation*}
for all $l=1,2,\dots,n$, so that $\frac{\Vert D^mf(0)(z_0^m)\Vert_\infty}{m!}\le1$. Moreover, we have
\begin{eqnarray}
  \frac{|D^{kp+m-1}g(0)(z_0^{kp+m-1})|}{(kp+m-1)!}\frac{|z_j|}{\Vert z\Vert_t}
  &\le& 1-\left(\frac{|D^{m-1}g(0)(z_0^{m-1})|}{(m-1)!}\frac{|z_j|}{\Vert z\Vert_t}\right)^2 \label{liu101}
\end{eqnarray}
for all $k\in\mathbb{N}_0$ and $z_0\in\partial B_{\ell_t^n}$.

On the other hand, if  $z_0\in\partial B_{\ell_t^n}$, then, for all $l=1,2,\dots,n$, we also have
\begin{equation}
\frac{|D^{k}f_l(0)(z_0^{k})|}{k!}=\frac{|D^{k-1}g(0)(z_0^{k-1})|}{(k-1)!}\frac{|z_l|}{\Vert z\Vert}_t\,
\mbox{~for~ all~}\, k\in\mathbb{N}.
\label{liu102}
\end{equation}

Combining (\ref{liu101}) and (\ref{liu102}), we find that
$$
\frac{|D^{kp+m}f_l(0)(z_0^{kp+m})|}{(kp+m)!}\le 1-\left(\frac{| D^{m}f_j(0)(z_0^{m})|}{m!}\right)^2,
$$
where $z_0\in\partial B_{\ell_t^n}$, $l=1,2,\dots,n$ and $k=1,2,\dots,n;$ that is,
$$
\frac{\Vert D^{kp+m}f(0)(z_0^{kp+m})\Vert_\infty}{(kp+m)!}\le 1-\left(\frac{\Vert D^{m}f(0)(z_0^{m})\Vert_\infty}{m!}\right)^2
$$
holds for $z_0\in\partial B_{\ell_t^n}$ and all $k\ge0$.

Because $z=z_0\Vert z\Vert$, by routine calculations, we obtain that
\beq
\sum_{k=1}^{\infty} \frac{\left \Vert D^{2kp+m}f(0)(z^{2kp+m})\right \Vert_\infty}{(2kp+m)!}
&=&\sum_{k=1}^{\infty} \frac{\left \Vert D^{2kp+m}f(0)(z_0^{2kp+m})\right \Vert_\infty}{(2kp+m)!}\Vert z\Vert_t^{2kp+m}\nonumber\\
&\le&\left[1-\left(\frac{\Vert  D^{m}f(0)(z_0^{m})\Vert_\infty}{m!}\right)^2\right]\sum_{k=1}^{\infty}\Vert z\Vert_t^{2kp+m}\nonumber\\
&=&\frac{\Vert z\Vert_t ^{2p-m}}{1-\Vert z\Vert_t ^{2p}}\left [ \Vert z\Vert_t ^{2m}-\left(\frac{\Vert  D^{m}f(0)(z^{m})\Vert_\infty}{m!}\right)^2\right ]\nonumber
\eeq
for all $\Vert z\Vert_t=r\in[0,1)$. This completes the proof.\hfill $\Box$

\vspace{8pt}

For  holomorphic mappings from $B_{\ell_t^n}$, $1\leq t\leq \infty$, to $\overline{\mathbb{D}}$, using the method of proof as in Lemma \ref{LLP-lem2}, we can easily derive the following, so we omit its proof.

\blem \label{LLP-lem1}
Suppose that $m\in\mathbb{N}_0$, $p\in\mathbb{N}$, $0\le m\le p$, and $1\leq t\leq \infty$. If $f \in H(B_{\ell_t^n},\overline{\mathbb{D}})$ and
$$f(z)=\frac{D^mf(0)(z^m)}{m!}+\sum_{k=1}^{\infty}\frac{D^{kp+m}f(0)(z^{kp+m})}{(kp+m)!},
$$
then
\begin{equation*}
\sum_{k=1}^{\infty} \frac{\left |D^{2kp+m}f(0)(z^{2kp+m})\right |}{(2kp+m)!}\le\frac{\Vert z\Vert_t ^{2p-m}}{1-\Vert z\Vert_t ^{2p}}\left [ \Vert z\Vert_t ^{2m}-\left(\frac{|D^m f(0)(z^m)|}{m!}\right)^2\right ]
\end{equation*}
holds for all $\Vert z\Vert_t=r\in[0,1)$.
\elem

\section{Bohr inequality for holomorphic mappings with lacunary series}\label{LLP-sec3}
In this section, we will use the Fr\'echet derivative to extend the Bohr inequality to higher dimensional space.

\subsection{Extension of Theorem B to the holomorphic mappings from $\mathbb{B}_X$ to $\overline{\mathbb{D}}$.}

\bthm\label{LLP-th01}
 Suppose that $f\in H(\mathbb{B}_X,\overline{\mathbb{D}})$, $f(z)=\sum_{k=0}^{\infty}\frac{D^{k}f(0)(z^{k})}{k!}$ and $p>0$. Then
\begin{equation}
|f(0)|^p+\sum_{k=1}^{\infty} \frac{\left |D^{k}f(0)(z^{k})\right |}{k!}+\left(\frac{1}{1+|f(0)|}+\frac{\Vert z\Vert}{1-\Vert z\Vert}\right)\sum_{k=1}^{\infty}\left(\frac{|D^{k} f(0)(z^{k})|}{k!}\right)^2\le1
\label{liu113b}
\end{equation}
holds for all $\Vert z\Vert=r\le r_p$, where $r_p=\frac{1-|f(0)|^p}{2-|f(0)|^2-|f(0)|^p}$. The number $r_p$ is best possible.

Moreover, when $p=1$, \eqref{liu113b} holds for $\Vert z\le \frac{1}{2+|f(0)|}$, and for $p=2$, it holds for $\|z\|\leq 1/2$.
\ethm

\noindent
{\bf Proof.}\quad
Fix $z\in \mathbb{B}_X\backslash\{ 0 \}$ and set $z_0=\frac{z}{\Vert z \Vert}$. For $\lambda\in \mathbb{D}$, we define $h(\lambda)=f(\lambda z_0).$ Then $h\in H(\mathbb{D},\overline{\mathbb{D}})$  and
\begin{equation*}
h(\lambda)
=b_0+\sum_{k=1}^{\infty}b_k\lambda^{k},
\end{equation*}
where
$$b_0= f(0)=\frac{D^{0}f(0)(z_0^{0})}{0!}  ~\mbox{ and }~ b_k=\frac{D^{k}f(0)(z_0^{k})}{k!}~ (k\ge1).
$$
By Theorem B, we know that
\begin{equation*}
|b_0|^p+\sum_{k=1}^{\infty}|b_k|\,|\lambda|^k+\left(\frac{1}{1+|b_0|}+\frac{|\lambda|}{1-|\lambda|}\right)\sum_{k=1}^{\infty}|b_k|^2|\lambda|^{2k}\le |b_0|^p+(1-|b_0|^2)\frac{|\lambda|}{1-|\lambda|}.
\end{equation*}

Set $\lambda=\|z\|=r$ and note that $z=z_0\Vert z\Vert$ and $|b_0|=|f(0)|\leq 1.$ Then the last relation gives
\begin{eqnarray*}
  && \hspace{-1.2cm}|f(0)|^p+\sum_{k=1}^{\infty} \frac{\left |D^{k}f(0)(z^{k})\right |}{k!}+\left(\frac{1}{1+|f(0)|}+\frac{\Vert z\Vert}{1-\Vert z\Vert}\right)\sum_{k=1}^{\infty}\left(\frac{|D^{k} f(0)(z^{k})|}{k!}\right)^2\label{liu113c}\\
  & =& |f(0)|^p+\sum_{k=1}^{\infty} \frac{\left |D^{k}f(0)(z_0^{k})\right |}{k!}\Vert z\Vert^k+\left(\frac{1}{1+|f(0)|}+\frac{\Vert z\Vert}{1-\Vert z\Vert}\right)\sum_{k=1}^{\infty}\left(\frac{|D^{k} f(0)(z_0^{k})|}{k!}\right)^2\Vert z\Vert^{2k}\nonumber\\
 &= & |b_0|^p+ \sum_{k=1}^{\infty}|b_k||\lambda|^k+\left(\frac{1}{1+|b_0|}+\frac{|\lambda|}{1-|\lambda|}\right)\sum_{k=1}^{\infty}|b_k|^2|\lambda|^{2k}\\
 &\leq &|f(0)|^p+(1-|f(0)|^2)\frac{r}{1-r},
\end{eqnarray*}
which is less than or equal to $1$ provided $\|z\|=r\leq r_p$. Thus, (\ref{liu113b}) holds for $\|z\|=r\leq r_p$, where $r_p=\frac{1-|f(0)|^p}{2-|f(0)|^2-|f(0)|^p}$.

Finally, we prove that inequality \eqref{liu113b} does not hold true for $x\in r_0\mathbb{B}_X$, where $r_0\in (r_p,1)$. We know that there exists a $c\in (0,1)$ and $v\in \partial \mathbb{B}_X$ such that $cr_0>r_p$ and
$$
c\sup \{||x||:x\in \partial \mathbb{B}_X\}<||v||.
$$
Now, we consider a function $f_1$ on $\mathbb{B}_X$ defined by
$$
f_1(x)=L_1\bigg(\cfrac{c\psi_v(x)}{||v||}\bigg),
$$
where $L_1(z)=(a-z)/(1-az),\,z\in \mathbb{D}$, $a\in [0,1)$, $\psi_v$ is a bounded linear functional on $X$ with $\psi_v(v)=||v||$ and $||\psi_v||=1$.
Choose $x=r_0v$ and we obtain that
\begin{eqnarray*}
  && \hspace{-2.5cm} |f_1(0)|^p+\sum_{k=1}^{\infty} \frac{\left |D^{k}f_1(0)(x^{k})\right |}{k!}+\left(\frac{1}{1+|f_1(0)|}+\frac{\Vert x\Vert}{1-\Vert x\Vert}\right)\sum_{k=1}^{\infty}\left(\frac{|D^{k} f_1(0)(x^{k})|}{k!}\right)^2\label{liu113c-a}\\
  & =& a^p+\sum_{k=1}^{\infty}(1-a^2)a^{k-1}c^kr_0^k+\left(\frac{1}{1+a}+\frac{cr_0}{1-cr_0}\right)\sum_{k=1}^{\infty}(1-a^2)^2a^{2k-2}c^{2k}r_0^{2k}\nonumber\\
 &=& a^p+(1-a^2)\frac{cr_0}{1-cr_0}>1
 \end{eqnarray*}
and the proof is complete.
\hfill $\Box$

\vspace{8pt}

In the following theorem, we determine the Bohr inequality for holomorphic functions, which fix the origin, with a lacunary series.

\bthm\label{LLP-th1}
Suppose that $m\in\mathbb{N}_0$, $p\in\mathbb{N},$ and $0\le m\le p.$ If $f \in H(\mathbb{B}_X,\overline{\mathbb{D}})$ and
$$f(z)=\sum_{k=1}^{\infty}\frac{D^{kp+m}f(0)(z^{kp+m})}{(kp+m)!},
$$
then
\begin{equation} \label{liu11}
\sum_{k=1}^{\infty} \frac{\left |D^{kp+m}f(0)(z^{kp+m})\right |}{(kp+m)!}
+\left(\frac{1}{\Vert z\Vert^{p+m}+\Lambda}+\frac{\Vert z\Vert ^{-m}}{1-\Vert z\Vert ^{p}}\right)\sum_{k=2}^{\infty}\left(\frac{|D^{kp+m} f(0)(z^{kp+m})|}{(kp+m)!}\right)^2\le1
\end{equation}
holds for all $0<\Vert z\Vert=r\le\tilde{r}_{p,m}$, where $\Lambda=\frac{|D^{p+m}f(0)(z^{p+m})|}{(p+m)!}$, and $\tilde{r}_{p,m}$ is the unique root in $(0, 1)$ of $G(r)=0,$ where
\begin{equation} \label{thm2-defn G(r)}
G(r)=5r^{2p+m}-2r^{p+m}+r^m+4r^p-4.
\end{equation}
For each $p$ and $m$, the number $\tilde{r}_{p,m}$  is best possible.
\ethm

\noindent
{\bf Proof.}\quad
As with the proof of Theorem \ref{LLP-th01}, we fix $z\in \mathbb{B}_X\backslash\{ 0 \}$ and let $z_0=\frac{z}{\Vert z \Vert}$. For $\lambda\in \mathbb{D}$,  define $h(\lambda)=f(\lambda z_0).$  Then $h\in H(\mathbb{D},\overline{\mathbb{D}})$ and
\begin{equation*}
h(\lambda)=\lambda^m\sum_{k=1}^{\infty}\frac{D^{kp+m}f(0)(z_0^{kp+m})}{(kp+m)!}(\lambda^{p})^k=\lambda^m\sum_{k=1}^{\infty}b_k(\lambda^{p})^k=:\lambda^m\varphi(\lambda^p),
\end{equation*}
where $b_k=\frac{D^{kp+m}f(0)(z_0^{kp+m})}{(kp+m)!}$,  $\varphi \in H(\mathbb{D},\overline{\mathbb{D}})$ with
$$
\varphi(\lambda)=\sum_{k=1}^{\infty}b_k \lambda^k=\lambda\sum_{k=1}^{\infty}b_k(\lambda)^{k-1}=\lambda \sum_{k=0}^{\infty}B_k \lambda^{k},
$$
and $B_k=b_{k+1}$. Clearly, $\varphi(\lambda^p)=\sum_{k=0}^{\infty}B_k(\lambda^{p})^{k}\in H(\mathbb{D},\overline{\mathbb{D}})$.
Then, according to Theorem B, we have
\begin{equation*}
\sum_{k=0}^{\infty}|B_k|\,|\lambda^p|^k+\left(\frac{1}{1+|B_0|}+\frac{|\lambda|^p}{1-|\lambda|^p} \right)\sum_{k=1}^{\infty}|B_k|^2|\lambda^p|^{2k}\le|B_0|+(1-|B_0|^2)\frac{|\lambda|^p}{1-|\lambda|^p},
\end{equation*}
which implies that
\begin{eqnarray*}
&&\sum_{k=0}^{\infty}\frac{|D^{(k+1)p+m}f(0)(z_0^{(k+1)p+m})|}{((k+1)p+m)!}|\lambda|^{pk}+\left(\frac{1}{1+\frac{|D^{p+m}f(0)z_0^{p+m}|}{(p+m)!}}+\frac{|\lambda|^p}{1-|\lambda|^p}\right)\times\nonumber\\
&& \hspace{2.5cm}\sum_{k=1}^{\infty}\left(\frac{|D^{(k+1)p+m}f(0)(z_0^{(k+1)p+m})|}{((k+1)p+m)!}\right)^2|\lambda|^{2kp}\le |B_0|+(1-|B_0|^2)\frac{|\lambda|^p}{1-|\lambda|^p}.
\end{eqnarray*}
Multiplying by $|\lambda|^{p+m}$ on both sides of the above inequality yields
\begin{eqnarray*}
&&\hspace{-.7cm}\sum_{k=0}^{\infty}\frac{|D^{(k+1)p+m}f(0)(z_0^{(k+1)p+m})|}{((k+1)p+m)!}|\lambda|^{(k+1)p+m}+\left(\frac{|\lambda|^{2p+2m}}{|\lambda|^{p+m}+\frac{|D^{p+m}f(0)z_0^{p+m}|}{(p+m)!}|\lambda|^{p+m}}+\frac{|\lambda|^{2p+m}}{1-|\lambda|^p}\right)\times\nonumber\\
&&\sum_{k=1}^{\infty}\left(\frac{|D^{(k+1)p+m}f(0)(z_0^{(k+1)p+m})|}{((k+1)p+m)!}\right)^2|\lambda|^{2kp}\le|B_0|\, |\lambda|^{p+m}+(1-|B_0|^2)\frac{|\lambda|^{2p+m}}{1-|\lambda|^p}.
\end{eqnarray*}
Taking $|\lambda|=\Vert z\Vert$, since $z=z_0\Vert z\Vert$, shows that
\begin{eqnarray*}
&&\hspace{-.5cm}\sum_{k=0}^{\infty} \frac{|D^{(k+1)p+m}f(0)(z^{(k+1)p+m})|}{((k+1)p+m)!}+\left(\frac{1}{\Vert z\Vert^{p+m}+\frac{|D^{p+m}f(0)z^{p+m}|}{(p+m)!}}+\frac{\Vert z\Vert^{-m}}{1-\Vert z\Vert^p}\right)\times\nonumber\\
&&\hspace{1.5cm}\sum_{k=1}^{\infty}\left(\frac{|D^{(k+1)p+m} f(0)(z^{(k+1)p+m})|}{((k+1)p+m)!}\right)^2\le|B_0|\, \Vert z\Vert^{p+m}+(1-|B_0|^2)\frac{\Vert z\Vert^{2p+m}}{1-\Vert z\Vert^p}.
\end{eqnarray*}
It is already obtained in \cite[Theorem 3.2]{LLP23} that
$$
|B_0|\, \Vert z\Vert^{p+m}+(1-|B_0|^2)\frac{\Vert z\Vert^{2p+m}}{1-\Vert z\Vert^p}\leq 1
$$
for all $0<\Vert z\Vert=r\le\tilde{r}_{p,m}$ and $\tilde{r}_{p,m}$ is the unique root in $(0, 1)$ of $G(r)=0,$ where $G(r)$ is given by \eqref{thm2-defn G(r)}.

Finally, we prove that inequality \eqref{liu11} does not hold true for $x\in r_0\mathbb{B}_X$, where $r_0\in (\tilde{r}_{p,m},1)$. We know that there exists a $c\in (0,1)$ and $v\in \partial \mathbb{B}_X$ such that $cr_0>\tilde{r}_{p,m}$ and
$$
c\sup \{||x||:x\in \partial \mathbb{B}_X\}<||v||.
$$

Now, we consider a function $f$ on $\mathbb{B}_X$ defined by
$$
f(x)=L_2\bigg(\cfrac{c\psi_v(x)}{||v||}\bigg),
$$
where
$$L_2(z)=z^{p+m}\left (\frac{a-z^p}{1-az^p}\right),\,z\in \mathbb{D} \mbox{ and } a\in [0,1),$$
$\psi_v$ is a bounded linear functional on $X$ with $\psi_v(v)=||v||$ and $||\psi_v||=1$.
Choosing $x=r_0v$, we get
\begin{eqnarray*}
\Lambda=\frac{\left |D^{p+m}f(0)(x^{p+m})\right |}{(p+m)!}&=a(cr_0)^{p+m}
\end{eqnarray*}
and 
$$
\frac{|D^{kp+m} f(0)(x^{kp+m})|}{(kp+m)!}=a^{k-2}(1-a^2)(cr_0)^{kp+m}.
$$
Thus, we have
\begin{eqnarray}
&&\hspace{-1.5cm} \sum_{k=1}^{\infty} \frac{\left |D^{kp+m}f(0)(x^{kp+m})\right |}{(kp+m)!}
+\left(\frac{1}{\Vert x\Vert^{p+m}+\Lambda}+\frac{\Vert x\Vert ^{-m}}{1-\Vert x\Vert ^{p}}\right)
  \sum_{k=2}^{\infty}\left(\frac{|D^{kp+m} f(0)(x^{kp+m})|}{(kp+m)!}\right)^2\nonumber\\
  &=&a(cr_0)^{m+p}+(1-a^2)\frac{(cr_0)^{2p+m}}{1-a(cr_0)^p}+\frac{(1-a^2)(1-a)(cr_0)^{3p+m}}{(1-(cr_0)^p)(1-a(cr_0)^p)}\nonumber\\
  &=&a(cr_0)^{p+m}+(1-a^2)\frac{(cr_0)^{2p+m}}{1-(cr_0)^p}. \label{LLL-8.1}
\end{eqnarray}
In the proof of \cite[Theorem 3.2]{LLP23}, it has been proved that
$$
a(x)^{p+m}+(1-a^2)\frac{(x)^{2p+m}}{1-(x)^p}>1,
$$
for  $x>\tilde r_{p,m}$.
This completes the proof of the sharpness of the constant $\tilde r_{p,m}$.
\hfill $\Box$

\vspace{8pt}
If $m=0,$ then \eqref{thm2-defn G(r)} reduces to $(r^p+1)(5r^p-3)=0$ and hence, $\tilde{r}_{p,0}=\sqrt[p]{3/5}.$ In particular, the case $m=0$ and $p=1$ in Theorem \ref{LLP-th1}, gives the following corollary.

\bcor\label{LLP-cor1}
Let $f \in H(\mathbb{B}_X,\overline{\mathbb{D}})$ and $f(z)=\sum_{k=1}^{\infty}\frac{D^{k}f(0)(z^{k})}{k!}$. Then
\begin{equation*}
\sum_{k=1}^{\infty} \frac{|D^{k}f(0)(z^{k})|}{k!}+\left(\frac{1}{\Vert z\Vert+|Df(0)(z)|}+\frac{1}{1-\Vert z\Vert}\right)\sum_{k=2}^{\infty}\left(\frac{|D^{k} f(0)(z^{k})|}{k!}\right)^2\le1,
\end{equation*}
holds for all $0<\Vert z\Vert=r\le 3/5$. The constant $3/5$ is best possible.
\ecor

\subsection{Bohr type inequality for holomorphic functions with lacunary series from $B_{\ell_t^n}$ to $\mathbb{D}^n$}
The proof of the following theorem uses the method of proof of \cite[Theorem 3.5]{LLP23}, which will be also used to prove Theorems \ref{LLP-th5} and \ref{LLP-th3}, respectively.

\bthm\label{LLP-th4}
Suppose that $m,\, p\in\mathbb{N}$, $1\le m\le p$, and $1\leq t \leq \infty$. Let $f \in H(B_{\ell_t^n},\overline{\mathbb{D}}^n)$ and
$$f(z)=zg(z)=\frac{D^mf(0)(z^m)}{m!}+\sum_{k=1}^{\infty}\frac{D^{kp+m}f(0)(z^{kp+m})}{(kp+m)!} ,
$$
where  $g\in H(B_{\ell_t^n},\mathbb{C})$. Then
\begin{equation}
	\sum_{k=1}^{\infty} \frac{\left \Vert D^{(2k-1)p+m}f(0)(z^{(2k-1)p+m})\right \Vert_\infty}{((2k-1)p+m)!}+\frac{\Vert z\Vert_t ^{p-m}}{1-\Vert z\Vert_t ^{2p}}\sum_{k=0}^{\infty}\left(\frac{\Vert D^{kp+m} f(0)(z^{kp+m})\Vert_\infty}{(kp+m)!}\right)^2\le1
	\label{liu110}
\end{equation}
holds for all $\Vert z\Vert_t=r\le r_{p,m}$, where $r_{p,m}$ is the unique root in $(0, 1)$ of the equation $r^{p+m}+r^{2p}-1=0$. For each $p$ and $m$, the number $r_{p,m}$  is best possible.
\ethm

\noindent
{\bf Proof.}\quad
Fix $z\in B_{\ell_t^n}\backslash\{ 0 \}$, and set $z_0=\frac{z}{\Vert z \Vert_t}.$ Then $z_0\in\partial B_{\ell_t^n}$. Let $j$ be such that $|z_j|=\Vert z\Vert_\infty=\max_{1\le l\le n}\{ |z_l| \}$.
Because $f(z)=zg(z)$, simple calculations yield
\begin{equation*}
	\frac{D^kf_j(0)(z_0^k)}{k!}=\frac{D^{k-1}g(0)(z_0^{k-1})}{(k-1)!}\frac{z_j}{\Vert z\Vert_t} \ \mbox{ for all $k\in\mathbb{N}$.}
\end{equation*}
Let $h_j(\lambda)=f_j(\lambda z_0)$ for $\lambda\in \mathbb{D}.$ Then $h_j\in H(\mathbb{D},\overline{\mathbb{D}})$ and
\begin{eqnarray*}
	h_j(\lambda)&=&\lambda^m\left[\frac{D^mf_j(0)(z_0^m)}{m!}+\sum_{k=1}^{\infty}\frac{D^{kp+m}f_j(0)(z_0^{kp+m})}{(kp+m)!}\lambda^{kp}\right]\\
	&=&\lambda^m\left[\frac{D^{m-1}g(0)(z_0^{m-1})}{(m-1)!}\frac{z_j}{\Vert z\Vert_t}+\sum_{k=1}^{\infty}\frac{D^{kp+m-1}g(0)(z_0^{kp+m-1})}{(kp+m-1)!}\frac{z_j}{\Vert z\Vert_t}\lambda^{kp}\right]\\
	&=&b_0\lambda^m+\sum_{k=1}^{\infty}b_k\lambda^{kp+m},
\end{eqnarray*}
where
$$b_0=\frac{D^{m-1}g(0)(z_0^{m-1})}{(m-1)!}\frac{z_j}{\Vert z\Vert_t} \ \mbox{ and }~ b_k=\frac{D^{kp+m-1}g(0)(z_0^{kp+m-1})}{(kp+m-1)!}\frac{z_j}{\Vert z\Vert_t} \ \mbox{ for $k\geq 1$.}
$$

By equation \eqref{LemD*} in Lemma D, and
$$\frac{|\lambda|^{2p}}{1-|\lambda|^{2p}}\sum_{k=0}^{\infty}|b_k|^2|\lambda|^{(2k-1)p}=\frac{|\lambda|^{p}}{1-|\lambda|^{2p}}\sum_{k=0}^{\infty}|b_k|^2|\lambda|^{2kp},
$$we have
\begin{equation}
	\sum_{k=1}^{\infty}|b_{2k-1}|\,|\lambda|^{(2k-1)p}+\frac{|\lambda|^p}{1-|\lambda|^{2p}}\sum_{k=0}^{\infty}|b_k|^2|\lambda|^{2kp}
	\le\frac{|\lambda|^p}{1-|\lambda|^{2p}}.
	\label{liu111}
\end{equation}
As $|z_j|=\Vert z\Vert_\infty=\max_{1\le l\le n}\{|z_l|\}$, substituting the expression of $b_k \ (k\ge0)$ into (\ref{liu111}), we have
\beq\label{liu112}
&&\sum_{k=1}^{\infty}\frac{|D^{(2k-1)p+m-1}g(0)(z_0^{(2k-1)p+m-1})|}{((2k-1)p+m-1)!}\frac{|z_j|}{\Vert z\Vert_t}|\lambda|^{(2k-1)p} \nonumber \\
&&\hspace{0.8cm}+\frac{|\lambda|^p}{1-|\lambda|^{2p}}\sum_{k=0}^{\infty}\left(\frac{|D^{kp+m-1}g(0)(z_0^{kp+m-1})|}{(kp+m-1)!}\frac{|z_j|}{\Vert z\Vert_t}\right)^2|\lambda|^{2kp}
\le\frac{|\lambda|^p}{1-|\lambda|^{2p}}
\eeq
for all $k\in\mathbb{N}_0$ and $z_0\in\partial B_{\ell_t^n}$.

Moreover, if $z_0\in\partial B_{\ell_t^n}$, we also have
\begin{equation}
	\frac{|D^{k}f_l(0)(z_0^{k})|}{k!}=\frac{|D^{k-1}g(0)(z_0^{k-1})|}{(k-1)!}\frac{|z_l|}{\Vert z\Vert_t}\le\frac{|D^{k-1}g(0)(z_0^{k-1})|}{(k-1)!}\frac{|z_j|}{\Vert z\Vert_t},\quad l=1,2,\dots,n.
	\label{liu113}
\end{equation}
Combining (\ref{liu112}) and (\ref{liu113}), for $z_0\in\partial B_{\ell_t^n}$, we have
\begin{eqnarray}\label{liu113a}
	&&\hspace{-2cm} \sum_{k=1}^{\infty}\frac{|D^{(2k-1)p+m}f_{l_k}(0)(z_0^{(2k-1)p+m})|}{((2k-1)p+m)!}|\lambda|^{(2k-1)p}\\
	&+&\frac{|\lambda|^p}{1-|\lambda|^{2p}}\sum_{k=0}^{\infty}\left(\frac{|D^{kp+m}f_{j_k}(0)(z_0^{kp+m})|}{(kp+m)!}\right)^2|\lambda|^{2kp}
	\le\frac{|\lambda|^p}{1-|\lambda|^{2p}}\nonumber
\end{eqnarray}
for all $k\in\mathbb{N}_0$ and where $l_k=1,2,\dots,n$ and $j_k=1,2,\dots,n$.

Multiplying  both sides of (\ref{liu113a}) by $|\lambda|^m$ and setting $|\lambda|=\Vert z\Vert_t$ so that $z=z_0\Vert z\Vert_t$, we have
\begin{equation*}
	\sum_{k=1}^{\infty}\frac{|D^{(2k-1)p+m}f_{l_k}(0)(z^{(2k-1)p+m})|}{((2k-1)p+m)!}+\frac{\Vert z\Vert_t^{p-m}}{1-\Vert z\Vert_t^{2p}}\sum_{k=0}^{\infty}\left(\frac{|D^{kp+m}f_{j_k}(0)(z^{kp+m})|}{(kp+m)!}\right)^2\le\frac{\Vert z\Vert_t^{p+m}}{1-\Vert z\Vert_t^{2p}}
\end{equation*}
for $z\in {\overline {B_{\ell_t^n}}}$, and $l_k=1,2,\dots,n$ and $j_k=1,2,\dots,n$. In other words,
\begin{equation*}
	\sum_{k=1}^{\infty} \frac{\left \Vert D^{(2k-1)p+m}f(0)(z^{(2k-1)p+m})\right \Vert_\infty}{((2k-1)p+m)!}+\frac{\Vert z\Vert_t ^{p-m}}{1-\Vert z\Vert_t ^{2p}}\sum_{k=0}^{\infty}\left(\frac{\Vert D^{kp+m} f(0)(z^{kp+m})\Vert_\infty}{(kp+m)!}\right)^2\le\frac{\Vert z\Vert_t^{p+m}}{1-\Vert z\Vert_t^{2p}}
\end{equation*}
which is less than or equal to $1$ provided  $r^{p+m}+r^{2p}-1\le0$, where $r=\Vert z\Vert_t$. The desired conclusion \eqref{liu110} follows for $\Vert z\Vert_t=r\leq r_{p,m}$, where $r_{p,m}$ is as in the statement.

To prove the sharpness, we just consider the functions $g$ and $f=zg \in H(B_{\ell_t^n}, \overline{\mathbb{D}}^n)$ given by
$$
g(z) 
=z_1^{m-1}\frac{a-z_1^p}{1-az_1^p} \ \mbox{ and } \ f(z)=\left(z_1^{m}\frac{a-z_1^p}{1-az_1^p}, z_2z_1^{m-1}\frac{a-z_1^p}{1-az_1^p},\dots,z_nz_1^{m-1}\frac{a-z_1^p}{1-az_1^p} \right)',
$$
where $z=(z_1,z_2,\dots,z_n)'$ and $a\in[0,1)$. In this case, let $z=(z_1,0,\dots,0)'$, which implies that $\Vert z\Vert_t=|z_1|=r$, and according to the definition of Fr\'echet derivative, we have
\begin{equation*}
	Df(0)(z)=\left(\frac{\partial f_l(0)}{\partial z_j}\right)_{1\le l,j\le n}(z_1,z_2,\dots,z_n)'.
\end{equation*}
Because $z=(z_1,0,\dots,0)'$, we have $Df(0)(z)=\left(\frac{\partial f_1(0)}{\partial z_1}z_1,0,\dots,0\right)$, and therefore,
$$\Vert Df(0)(z)\Vert=\left|\frac{\partial f_1(0)}{\partial z_1}z_1\right|.
$$
With the help of the proof of Theorem \cite[Theorem 3.5]{LLP23}, we obtain
$$\frac{\Vert D^{kp+m}f(0)(z^{kp+m})\Vert_\infty}{(kp+m)!}=\left|\frac{\partial^{kp+m}f_1(0)}{\partial z_1^{kp+m}}\frac{z_1^{kp+m}}{(kp+m)!}\right| \quad \mbox{for $k\ge0$.}
$$
Therefore, we compute that
\begin{eqnarray*}
	&&\hspace{-.5cm}\sum_{k=1}^{\infty} \frac{\left \Vert D^{(2k-1)p+m}f(0)(z^{(2k-1)p+m})\right \Vert_\infty}{((2k-1)p+m)!}+\frac{\Vert z\Vert_t ^{p-m}}{1-\Vert z\Vert_t ^{2p}}\sum_{k=0}^{\infty}\left(\frac{\Vert D^{kp+m} f(0)(z^{kp+m})\Vert_\infty}{(kp+m)!}\right)^2\nonumber\\
	&&=\sum_{k=1}^{\infty}(1-a^2)a^{2(k-1)}r^{(2k-1)p+m}+\frac{r^{p-m}}{1-r^{2p}}a^2r^{2m}+\frac{r ^{p-m}}{1-r ^{2p}}\sum_{k=1}^{\infty}(1-a^2)^2a^{2(k-1)}r^{2kp+2m}\nonumber\\
	&&=(1-a^2)r^{p+m}\left[\frac{ (1-r^{2p}) -(1-a^2r^{2p}) +(1-a^2)r^{2p}}{(1-r^{2p})(1-a^2r^{2p})}\right] + \frac{r^{p+m}}{1-r^{2p}}=\frac{r^{p+m}}{1-r^{2p}},
\end{eqnarray*}
which is clearly bigger than $1$ whenever $r>r_{p,m}$.
This completes the proof. \hfill $\Box$

\vspace{8pt}

Using the method of proof Theorems \ref{LLP-th01} and \ref{LLP-th4}, we may verify the proof of the following, and so we omit its proof.

\bthm\label{LLP-th2}
Suppose that $m\in\mathbb{N}_0$, $p\in\mathbb{N}$,  $0\le m\le p$, and $1\leq t\leq \infty$. If $f \in H(B_{\ell_t^n},\overline{\mathbb{D}})$ and
$$
f(z)=\frac{D^mf(0)(z^m)}{m!}+\sum_{k=1}^{\infty}\frac{D^{kp+m}f(0)(z^{kp+m})}{(kp+m)!},
$$
then
\begin{equation}
	\sum_{k=1}^{\infty} \frac{\left |D^{(2k-1)p+m}f(0)(z^{(2k-1)p+m})\right |}{((2k-1)p+m)!}+\frac{\Vert z\Vert_t ^{p-m}}{1-\Vert z\Vert_t ^{2p}}\sum_{k=0}^{\infty}\left(\frac{|D^{kp+m} f(0)(z^{kp+m})|}{(kp+m)!}\right)^2\le1
	\label{liu14}
\end{equation}
holds for all $\Vert z\Vert_t=r\le r_{p,m}$, where $r_{p,m}$ is the same as in Theorem $\ref{LLP-th4}$.
The constant $r_{p,m}$ is best possible for each $p$ and $m$.
\ethm

\section{Bohr inequality for holomorphic mappings with alternating series}\label{LLP-sec4}
In this section, we will use the Fr\'echet derivative to extend the Bohr inequality with alternating series to the higher-dimensional space, and obtain
the higher-dimensional generalizations of Theorem C.

\bthm\label{LLP-th5}
Suppose that $m,\, p\in\mathbb{N}$, $p$ is odd, $1\le m\le p$, and $1\leq t\leq \infty$. If $f \in H(B_{\ell_t^n},\overline {\mathbb{D}}^n)$ and
$$
f(z)=zg(z)=\frac{D^mf(0)(z^m)}{m!}+\sum_{k=1}^{\infty}\frac{D^{kp+m}f(0)(z^{kp+m})}{(kp+m)!},
$$
where $g\in H(B_{\ell_t^n},\mathbb{C}),$ then
\begin{eqnarray}\label{liu117}
	&&\hspace{-2cm}\Bigg| \sum_{k=1}^{\infty}(-1)^{kp+m} \frac{\left \Vert D^{kp+m}f(0)(z^{kp+m})\right \Vert_\infty}{(kp+m)!}\\
	&&+ \,(-1)^{m+p}\frac{\Vert z\Vert_t ^{p-m}}{1-\Vert z\Vert_t ^{2p}}\sum_{k=0}^{\infty}\left(\frac{\Vert D^{kp+m} f(0)(z^{kp+m})\Vert_\infty}{(kp+m)!}\right)^2\Bigg| \leq 1\nonumber
\end{eqnarray}
holds for all $\Vert z\Vert_t=r\le r_{p,m}$, where $r_{p,m}$ is the same as in Theorem $\ref{LLP-th4}$.
The constant $r_{p,m}$ is best possible for each $p$ and $m$.
\ethm

\noindent
{\bf Proof.}\quad
(i) Assume first that $p$ is odd and $m$ is even. Then $(-1)^{m+p}=-1$, and
\begin{eqnarray}\label{liu118}
	&&\sum_{k=1}^{\infty}(-1)^{kp+m} \frac{\left \Vert D^{kp+m}f(0)(z^{kp+m})\right \Vert_\infty}{(kp+m)!}=\sum_{k=1}^{\infty} \frac{\left \Vert D^{2kp+m}f(0)(z^{2kp+m})\right \Vert_\infty}{(2kp+m)!}\\
	&& \hspace{7.3cm}-\sum_{k=1}^{\infty} \frac{\left \Vert D^{(2k-1)p+m}f(0)(z^{(2k-1)p+m})\right \Vert_\infty}{((2k-1)p+m)!}.\nonumber
\end{eqnarray}

Now, to find the lower bound, by Theorem \ref{LLP-th4} and (\ref{liu118}), we see that
\beq
&&\hspace{-1cm}\sum_{k=1}^{\infty}(-1)^{kp+m} \frac{\left \Vert D^{kp+m}f(0)(z^{kp+m})\right \Vert_\infty }{(kp+m)!}-\frac{\Vert z\Vert_t ^{p-m}}{1-\Vert z\Vert_t ^{2p}}\sum_{k=0}^{\infty}\left(\frac{\Vert D^{kp+m} f(0)(z^{kp+m})\Vert_\infty }{(kp+m)!}\right)^2 \nonumber\\
&\ge&-\sum_{k=1}^{\infty} \frac{\left \Vert D^{(2k-1)p+m}f(0)(z^{(2k-1)p+m})\right \Vert_\infty }{((2k-1)p+m)!}-\frac{\Vert z\Vert_t ^{p-m}}{1-\Vert z\Vert_t ^{2p}}\sum_{k=0}^{\infty}\left(\frac{\Vert D^{kp+m} f(0)(z^{kp+m})\Vert_\infty }{(kp+m)!}\right)^2\nonumber\\
&\ge&-1\nonumber
\eeq
holds for all $\Vert z\Vert_t=r\le r_{p,m}$, where $r_{p,m}$ is the unique root in $(0, 1)$ of the equation $r^{p+m}+r^{2p}-1=0$.

To find the upper bound, according to Lemma \ref{LLP-lem2} and the method of proof of Lemma \ref{LLP-lem2}, we see that $\frac{\Vert D^m f(0)(z_0^m)\Vert_\infty}{m!}\le1$ and
\beq
&&\sum_{k=1}^{\infty}(-1)^{kp+m} \frac{\left \Vert D^{kp+m}f(0)(z^{kp+m})\right \Vert_\infty }{(kp+m)!}-\frac{\Vert z\Vert_t ^{p-m}}{1-\Vert z\Vert_t ^{2p}}\sum_{k=0}^{\infty}\left(\frac{\Vert D^{kp+m} f(0)(z^{kp+m})\Vert_\infty }{(kp+m)!}\right)^2 \label{liu120}\\
&\le&\sum_{k=1}^{\infty} \frac{\left \Vert D^{2kp+m}f(0)(z^{2kp+m})\right \Vert_\infty }{(2kp+m)!}\le\frac{\Vert z\Vert_t ^{2p-m}}{1-\Vert z\Vert_t ^{2p}}\left [ \Vert z\Vert_t ^{2m}-\left(\frac{\Vert D^m f(0)(z^m)\Vert_\infty}{m!}\right)^2\right ]\nonumber\\
&=&\frac{\Vert z\Vert_t ^{2p+m}}{1-\Vert z\Vert_t ^{2p}}\left [ 1-\left(\frac{\Vert D^m f(0)(z_0^m)\Vert_\infty }{m!}\right)^2\right ]\le\frac{\Vert z\Vert_t ^{2p+m}}{1-\Vert z\Vert_t ^{2p}}\le\frac{r^{p+m}}{1-r ^{2p}}\nonumber,
\eeq
which is less than or equal to 1 whenever
$\frac{r^{p+m}}{1-r ^{2p}}\leq 1$, that is, whenever $r^{p+m}+r ^{2p}-1\le0$. Thus, combining with the value of the upper and lower bounds, (\ref{liu117}) holds for all $\Vert z\Vert_t=r\le r_{p,m}$, where $r_{p,m}$ is the unique root in $(0, 1)$ of the equation $r^{p+m}+r^{2p}-1=0$.

To prove the sharpness, we just consider the functions $g(z)=z_1^{m+p-1}$, and $f \in H(B_{\ell_t^n}, \overline{\mathbb{D}}^n)$ given by
$f(z)=zg(z)=(z_1^{m+p}, z_2z_1^{m+p-1},\dots,z_nz_1^{m+p-1})'$, where $z=(z_1,z_2,\dots,z_n)'$.

Let $z=(z_1,0,\dots,0)'$. Then $\Vert z\Vert_t=|z_1|=r$. In this case, just like the method of proof of Theorem \ref{LLP-th4}, we have $$
\frac{\left \Vert D^{kp+m}f(0)(z^{kp+m})\right \Vert_\infty}{(kp+m)!}=\left|\frac{\partial^{kp+m}f_1(0)}{\partial z_1^{kp+m}}\frac{z_1^{kp+m}}{(kp+m)!}\right| \ \mbox{for all $k\ge0$}.
$$
Since $f_1(z)=z_1^{m+p}$, we have, $\frac{\left \Vert D^{kp+m}f(0)(z^{kp+m})\right \Vert_\infty}{(kp+m)!}=0$ (when $k\neq 1$), and $m+p$ is odd. Therefore,
\begin{eqnarray*}
	&&\hspace{-.7cm}\left | \sum_{k=1}^{\infty}(-1)^{kp+m} \frac{\left \Vert D^{kp+m}f(0)(z^{kp+m})\right \Vert_\infty}{(kp+m)!}+(-1)^{m+p}\frac{\Vert z\Vert_t ^{p-m}}{1-\Vert z\Vert_t ^{2p}}\sum_{k=0}^{\infty}\left(\frac{\Vert D^{kp+m} f(0)(z^{kp+m})\Vert_\infty}{(kp+m)!}\right)^2\right |\nonumber\\
	&&=\left| -r^{m+p}-\frac{r^{p-m}}{1-r^{2p}}r^{2(p+m)}\right|=\frac{r^{p+m}}{1-r^{2p}},
\end{eqnarray*}
which shows that the left hand of (\ref{liu117}) is bigger than $1$ whenever $r>r_{p, m}$. This completes the proof of Case (i).
\vskip 2mm

(ii) Next, suppose that both $p$ and $m$ are odd. Then we see that $(-1)^{m+p}=1$ and
\begin{eqnarray*}
	&&\sum_{k=1}^{\infty}(-1)^{kp+m} \frac{\left \Vert D^{kp+m}f(0)(z^{kp+m})\right \Vert_\infty}{(kp+m)!}=-\sum_{k=1}^{\infty} \frac{\left \Vert D^{2kp+m}f(0)(z^{2kp+m})\right \Vert_\infty}{(2kp+m)!}\nonumber\\
	&& \hspace{7.3cm}+\sum_{k=1}^{\infty} \frac{\left \Vert D^{(2k-1)p+m}f(0)(z^{(2k-1)p+m})\right \Vert_\infty}{((2k-1)p+m)!}.
\end{eqnarray*}

Therefore, combined with Lemma \ref{LLP-lem2} and Theorem \ref{LLP-th4}, the rest of the proof is similar to the first part above, except that the corresponding formula has the opposite sign. This completes the proof. \hfill $\Box$

\vskip 2mm

In the following case, since $f((0,0,\dots,0)')=(0,0,\dots,0)' \;g((0,0,\dots,0)')=(0,0,\dots,0)'$, Corollary \ref{LLP-cor3} is an extension of Theorem B and Corollary \ref{LLP-cor1}.

\bcor\label{LLP-cor3}
Let $f \in H(B_{\ell_t^n},\overline{\mathbb{D}}^n)$ and $f(z)=zg(z)=\sum_{k=1}^{\infty}\frac{D^{k}f(0)(z^{k})}{k!}$, where $g\in H(B_{\ell_t^n},\mathbb{C})$, and $j$ satisfies $|z_j|=\Vert z\Vert_\infty=\max_{1\le l\le n}\{ |z_l| \}$. Then
\begin{equation}
	\sum_{k=1}^{\infty} \frac{\Vert D^{k}f(0)(z^{k})\Vert_\infty}{k!}+\left(\frac{1}{\Vert z\Vert_t+\Vert Df(0)(z)\Vert_\infty}+\frac{1}{1-\Vert z\Vert_t}\right)\sum_{k=2}^{\infty}\left(\frac{\Vert D^{k} f(0)(z^{k})\Vert_\infty}{k!}\right)^2\le1,
	\label{liu119}
\end{equation}
holds for all $0<\Vert z\Vert_t=r\le 3/5$. This result is sharp.
\ecor

\noindent
{\bf Proof.} 
Fix $z\in B_{\ell_t^n}\backslash\{ 0 \}$, and set $z_0=\frac{z}{\Vert z \Vert_t}\in\partial B_{\ell_t^n}$. Because $f(z)=zg(z)$, through simple calculations, we have
\begin{equation*}
	\frac{D^kf_l(0)(z_0^k)}{k!}=\frac{D^{k-1}g(0)(z_0^{k-1})}{(k-1)!}\frac{z_l}{\Vert z\Vert_t}
\end{equation*}
for all $l=1,2,\dots,n$ and $k\in\mathbb{N}_0$. Since $|z_j|=\Vert z\Vert_\infty=\max_{1\le l\le n}\{|z_l|\}$, it follows that
\begin{equation*}
	\frac{|D^{k}f_l(0)(z_0^{k})|}{k!}\le\frac{|D^{k-1}g(0)(z_0^{k-1})|}{(k-1)!}\frac{|z_j|}{\Vert z\Vert_t}=\frac{|D^{k}f_j(0)(z_0^{k})|}{k!}
\end{equation*}
for all $l=1,2,\dots,n$ and $k\in\mathbb{N}_0$. Therefore,
$$\frac{\Vert D^{k}f(0)(z_0^{k})\Vert_\infty}{k!}=\frac{|D^{k}f_j(0)(z_0^{k})|}{k!} \ \mbox{for all $k\in\mathbb{N}_0.$}
$$
Multiplying both sides of the above equation by $\Vert z^k\Vert_t$ gives
$$
\frac{\Vert D^{k}f(0)(z^{k})\Vert_\infty}{k!}=\frac{|D^{k}f_j(0)(z^{k})|}{k!}\ \mbox{for all $k\in\mathbb{N}_0$.}
$$
Therefore,
\begin{eqnarray*}
	&&\sum_{k=1}^{\infty} \frac{\Vert D^{k}f(0)(z^{k})\Vert_\infty}{k!}+\left(\frac{1}{\Vert z\Vert_t+\Vert Df(0)(z)\Vert_\infty}+\frac{1}{1-\Vert z\Vert_t}\right)\sum_{k=2}^{\infty}\left(\frac{\Vert D^{k} f(0)(z^{k})\Vert_\infty}{k!}\right)^2\nonumber\\
	&& \hspace{0.5cm}=\sum_{k=1}^{\infty} \frac{|D^{k}f_j(0)(z^{k})|}{k!}+\left(\frac{1}{\Vert z\Vert_t+|Df_j(0)(z)|}+\frac{1}{1-\Vert z\Vert_t}\right)\sum_{k=2}^{\infty}\left(\frac{| D^{k} f_j(0)(z^{k})|}{k!}\right)^2.
\end{eqnarray*}
Since $f_j\in H(B_{\ell_t^n},\overline{\mathbb{D}})$, by Corollary \ref{LLP-cor1}, we  find that $(\ref{liu119})$ holds for all $0<\Vert z\Vert_t=r\le 3/5$.

To prove the sharpness, we just consider the function $g(z)
=\frac{a-z_1}{1-az_1}$ and $f \in H(B_{\ell_t^n}, \overline{\mathbb{D}}^n)$ given by
$$f(z)=zg(z)= \left (z_1\frac{a-z_1}{1-az_1},z_2\frac{a-z_1}{1-az_1},\dots,z_n\frac{a-z_1}{1-az_1}\right )'.
$$
In this case, let $z=(z_1,0,\dots,0)'$. The rest of the proof is similar as in Theorem \ref{LLP-th4}.\hfill $\Box$

\bcor\label{LLP-cor4} 
Suppose that $m,\, p\in\mathbb{N}$, $p$ is odd, $1\le m\le p$, and $1\leq t\leq \infty$. Let $f  \in H(B_{\ell_t^n},\overline {\mathbb{D}}^n)$ be given by $$f(z)=zg(z)=\frac{D^mf(0)(z^m)}{m!}+\sum_{k=1}^{\infty}\frac{D^{kp+m}f(0)(z^{kp+m})}{(kp+m)!},
$$
where $g\in H(B_{\ell_t^n},\mathbb{C})$. Then
\begin{eqnarray*}
	\left|A_{f_m}(r)+(-1)^m\left(\frac{1}{\Vert z\Vert_t^{m}+\Gamma}+\frac{\Vert z\Vert_t ^{2p-m}}{1-\Vert z\Vert_t ^{2p}}\right)\sum_{k=0}^{\infty}\left(\frac{\Vert D^{kp+m} f(0)(z^{kp+m})\Vert_\infty}{(kp+m)!}\right)^2\right|\leq 1
\end{eqnarray*}
holds for $0<\|z\|_t = r\leq \tilde{R}_{p,m}$, where
$$
A_{f_m}(r)=\sum_{k=1}^{\infty}(-1)^{kp+m}\frac{\Vert D^{kp+m}f(0)(z^{kp+m})\Vert_\infty}{(kp+m)!},\quad \Gamma=\frac{\Vert D^m f(0)(z^m)\Vert_\infty}{m!},
$$
and $\tilde{R}_{p,m}$ is the unique root in $(0,1)$ of the equation $r^{2p+m}+2r^{2p}-1=0$. For each $p, m,$ the constant $\tilde{R}_{p,m}$ is best possible.
\ecor

\noindent
{\bf Proof.}\quad
First, use the methods of proofs of Lemma \ref{LLP-lem2} and Theorem \ref{LLP-th4} to obtain the upper bounds of odd and even terms respectively, and then use (\ref{liu21b}) of Lemma D and
combine the  methods of proofs Corollary \ref{LLP-cor3} and Theorem \ref{LLP-th5} to get quickly the proof of Corollary \ref{LLP-cor4}, so we skip the details.\hfill $\Box$

\br\label{HLP-re2}
When $f \in H(B_{\ell_t^n},\overline{\mathbb{D}})$, $m,\, p\in\mathbb{N}_0$, $p$ is odd, $0\le m\le p$, and $1\leq t\leq \infty$, then the same conclusion as Corollary \ref{LLP-cor4} can be obtained. In the case of $m=0$, we have $\Vert z\Vert_t=r\le 1/\sqrt[2p]{3}$, and thus this case can be regarded as an extension of Theorem 1.2 in \cite{ABS2017}.
\er
\vskip 2mm

Combining Lemma \ref{LLP-lem1} and Theorem \ref{LLP-th2}, and using the method of proof of Theorem \ref{LLP-th5}, we can easily obtain the proof of the following theorem.

\bthm\label{LLP-th3}
Suppose that $m\in\mathbb{N}_0$, $p\in\mathbb{N}$, $p$ is odd, $0\le m\le p$ and $1\leq t\leq \infty$. If $f  \in H(B_{\ell_t^n},\overline{\mathbb{D}})$ and
$$
f(z)=\frac{D^mf(0)(z^m)}{m!}+\sum_{k=1}^{\infty}\frac{D^{kp+m}f(0)(z^{kp+m})}{(kp+m)!},
$$
then
\begin{equation*}
	\left | \sum_{k=1}^{\infty}(-1)^{kp+m} \frac{\left |D^{kp+m}f(0)(z^{kp+m})\right |}{(kp+m)!}+(-1)^{m+p}\frac{\Vert z\Vert_t ^{p-m}}{1-\Vert z\Vert_t ^{2p}}\sum_{k=0}^{\infty}\left(\frac{|D^{kp+m} f(0)(z^{kp+m})|}{(kp+m)!}\right)^2\right | \leq 1
\end{equation*}
holds for all $\Vert z\Vert_t=r\le r_{p,m}$, where $r_{p,m}$ is the same as in Theorem \ref{LLP-th4}.
For each $p$ and $m$ the number $r_{p,m}$ is best possible.
\ethm

\br\label{HLP-re3}
When $n=1$ (that is $f\in H(\mathbb{D},\overline{\mathbb{D}})$), through comparison, it can be easily found that the results of Theorems \ref{LLP-th5} and \ref{LLP-th3} are the same as that of Theorem C.
\er

\subsection*{Acknowledgements}
 The work of the first author is supported by the institute postdoctoral fellowship of IIT Madras, India. 
 All authors contribute equally.

\subsection*{Conflict of Interests}
The authors declare that they have no conflict of interest, regarding the publication of this paper.

\subsection*{Data Availability Statement}
The authors declare that this research is purely theoretical and does not associate with any datas.

\end{document}